\documentclass[a4paper]{amsart}

\usepackage{amssymb,latexsym}
\usepackage{amsmath}
\usepackage{amsfonts}
\usepackage{amsthm}
\usepackage{amssymb}

\theoremstyle{plain}
\newtheorem{theorem}{Theorem}

\newtheorem{proposition}[theorem]{Proposition}

\theoremstyle{definition}

\newtheorem{example}{Example}
\newtheorem{remark}{Remark}

\title{Strong Eulerian triples}
\begin{document}

\date{}


\author[A. Dujella]{Andrej Dujella}
\address{
Department of Mathematics\\
University of Zagreb\\
Bijeni{\v c}ka cesta 30, 10000 Zagreb, Croatia
}
\email[A. Dujella]{duje@math.hr}

\author[I. Gusi\'c]{Ivica Gusi\'c}
\address{
Faculty of Chemical Engin. and Techn. \\
University of Zagreb \\
Maruli\'cev trg 19, 10000 Zagreb, Croatia
}
\email[I. Gusi\'c]{igusic@fkit.hr}

\author[V. Petri\v{c}evi\'c]{Vinko Petri\v{c}evi\'c}
\address{
Department of Mathematics\\
University of Zagreb\\
Bijeni{\v c}ka cesta 30, 10000 Zagreb, Croatia
}
\email[V. Petri\v{c}evi\'c]{vpetrice@math.hr}

\author[P. Tadi\'c]{Petra Tadi\'c}
\address{
Department of Economics and Tourism \\
Juraj Dobrila University of Pula \\
52100 Pula, Croatia
}
\email[P. Tadi\'c]{ptadic@unipu.hr}

\begin{abstract}
We prove that there exist infinitely many rationals $a$, $b$ and $c$ with the property that
$a^2-1$, $b^2-1$, $c^2-1$, $ab-1$, $ac-1$ and $bc-1$ are all perfect squares. This provides a solution
to a variant of the problem studied by Diophantus and Euler.
\end{abstract}


\maketitle

\section{Introduction}

Let $q$ be a non-zero rational. A set $\{a_1, a_2, \ldots, a_m\}$ of $m$ non-zero rationals
is called a rational $D(q)$-$m$-tuple if $a_i \cdot a_j + q$ is a perfect square for all
$1 \leq· i < j \leq m$. Diophantus found the first rational $D(1)$-quadruple
$\{1/16, 33/16, 17/4, 105/16\}$, while Euler found a rational $D(1)$-quintuple by
extending the integer $D(1)$-quadruple $\{1, 3, 8, 120\}$, found by Fermat,
with fifth rational number $777480/8288641$ (see \cite{Dickson,Heath}).
Recently, Stoll \cite{Stoll} proved that extension of Fermat's set to a rational $D(1)$-quintuple
is unique. The first example of a rational $D(1)$-sextuple, the set
$\{11/192, 35/192, 155/27, 512/27, 1235/48, 180873/16\}$, was found by Gibbs \cite{Gibbs1},
while Dujella, Kazalicki, Miki\'c and Szikszai \cite{DKMS} recently proved that there are
infinitely many rational $D(1)$-sextuples (see also \cite{Duje-Matija}).
It is not known whether there exist any rational $D(1)$-septuple. However, Gibbs \cite{Gibbs2}
found examples of ``almost'' septuples, namely, rational $D(1)$-quintuples which can be extended
to rational $D(1)$-sextuples on two different ways, so that only one condition is missing that these
seven numbers form a rational $D(1)$-septuple (they form $D(1)$-septuples over corresponding
quadratic fields). One such quintuple is
$\{243/560, 1147/5040, 1100/63$, $7820/567, 95/112\}$ which can be extended to the sextuple with
$38269/6480$ or $196/45$. For an overview of results on $D(1)$-$m$-tuples and its generalizations
see \cite{Duje-Notices}.

It is known that for every rational $q$ there exist infinitely many rational $D(q)$-quadruples
(see \cite{Duje-graz2000}). In 2012, Dujella and Fuchs \cite{Duje-Fuchs-JNT} proved that for
infinitely many square-free integers $q$ there are infinitely many rational $D(q)$-quintuples,
by considering twists of the elliptic curve
$y^2 = x^3 + 86x^2 + 825x$ with positive rank.

Apart of the case $q=1$, the most studied case in the literature is $q=-1$.
The case $q=-1$ is closely related to another old problem investigated by Diophantus and Euler,
concerning the sets of integers or rationals with the property that the product of any two
of its distinct elements plus their sum is a perfect square. We call a set
$\{x_1,x_2,\ldots,x_m\}$ an Eulerian $m$-tuple if $x_ix_j+x_i+x_j$ is a perfect square for
all $1 \leq· i < j \leq m$. The equality $x_ix_j+x_i+x_j = (x_i+1)(x_j+1)-1$ gives explicit
connection between Eulerian $m$-tuples and $D(-1)$-$m$-tuples.
It is known that there does not exist a $D(-1)$-quintuple in integers and that there
are only finitely many such quadruples, and all of them have to contain the element $1$
(see \cite{DFF,Duje-Fuchs}). In particular, there does not exist an Eulerian quadruple in positive integers. On the other hand, there exist infinitely many rational $D(-1)$-quintuples, and hence infinitely
many Eulerian quintuples in rationals (see \cite{Duje-Euler1,Duje-Euler2}).

Note that in the definition of rational $D(q)$-$m$-tuples we excluded the requirement
that the product of an element with itself plus $q$ is a square. It is obvious that for $q=1$
such condition cannot be satisfied in integers. But for rationals there is no obvious reason why the sets
(called strong $D(1)$-$m$-tuples) which satisfy these stronger conditions would not exist.
Dujella and Petri\v{c}evi\'c \cite{Duje-Vinko} proved in
2008 that there exist infinitely many strong $D(1)$-triples, while no example of a strong
$D(1)$-quadruple is known.

In this paper, we study the existence of strong Eulerian triples, i.e. sets of three rationals
$\{x_1,x_2,x_3\}$ such that $x_1x_2+x_1+x_2$, $x_1x_3+x_1+x_3$, $x_2x_3+x_2+x_3$,
$x_1^2+2x_1$, $x_2^2+2x_2$ and $x_3^2+2x_3$ are all perfect squares.
Equivalently, by taking $a_i=x_i+1$, we may consider strong rational $D(-1)$-triples,
i.e. sets of three rationals $\{a_1,a_2,a_3\}$ such that $a_1a_2-1$, $a_1a_3-1$, $a_2a_3-1$,
$a_1^2-1$, $a_2^2-1$ and $a_3^2-1$ are all perfect squares.
It is clear that all elements of a strong rational $D(-1)$-triple has to have the same sign,
and that $\{a_1,a_2,a_3\}$ is a strong rational $D(-1)$-triple if and only if $\{-a_1,-a_2,-a_3\}$
has the same property. Thus, there is no loss of generality in assuming that
all elements of a strong rational $D(-1)$-triple are positive.
By connecting the problem with certain families of elliptic curves, we will show that
there exist infinitely many strong rational $D(-1)$-triples.
We find only eight strong rational $D(-1)$-triples that do not contain the number $1$
(see Example \ref{ex:1} and Remark \ref{rem3}).
Accordingly, our construction gives several infinite families of strong rational $D(-1)$-triples
which all contain the number $1$.
This means that the corresponding strong Eulerian triples contain the number $0$,
and all other elements are squares.
MacLeod \cite{MacLeod} found examples of rational Eulerian triples and quadruples
which all elements are squares. However, in our situation there is an additional requirement
that each element increased by $2$ is also a square.

The main result of this paper, which will be proved in Section \ref{sec:2}, is the
following theorem.

\begin{theorem} \label{ref:1}
There exist infinitely many strong Eulerian triples.
\end{theorem}

This is a more precise formulation of our result.

\begin{proposition} \label{ref:2}
Let $\{1,b\}$ be a strong rational $D(-1)$-pair.
Then there exist infinitely many strong rational $D(-1)$-triples of the form $\{1,b,c\}$.
\end{proposition}

\section{Constructions of infinite families of triples} \label{sec:2}

\begin{example}  \label{ex:1}
We start by searching experimentally for strong rational $D(-1)$-triples
with elements with relatively small numerators and denominators (smaller than $2.5\cdot 10^7$).
We found seven examples with all elements different from $1$:
$$ \{493/468, 1313/1088, 33137/32912\}, $$
$$ \{1517/1508, 42601/11849, 909745/757393\}, $$
$$ \{125/117, 689/400, 14353373/13130325\}, $$
$$ \{354005/22707, 193397/183315, 2084693/2074035\}, $$
$$ \{2833349/218660, 3484973/2619045, 3056365/3047653\}, $$
$$ \{2257/1105, 2873/2745, 3859145/862784\}, $$
$$ \{2257/1105, 115825/8177, 14307761/10303760\}, $$
and $23$ examples containing the number $1$:

\begin{eqnarray*}
& \{1, 5/4, 14645/484\}, \quad \{1, 689/400, 1025/64\}, \\
& \{1, 689/400, 969425/861184\}, \quad \{1, 689/400, 9047825/4857616\}, \\
& \{1, 2501/100, 59189/12100\}, \quad \{1, 2501/100, 3219749/2102500\}, \\
& \{1, 6625/1296, 3254641/435600\}, \quad \{1, 19825/17424, 46561/32400\}, \\
& \{1, 19825/17424, 50689/3600\}, \quad \{1, 17009/6400,  8530481/4494400\}, \\
& \{1, 26245/324,  26361205/18301284\}, \quad \{1, 28625/2704,  27060449/25603600\}, \\
& \{1, 60229/44100, 65125/39204\}, \quad \{1, 65125/39204,  2829205/30276\}, \\
& \{1, 168305/94864, 262145/1024\}, \quad \{1, 926021/96100,  13236725/7365796\}, \\
& \{1, 1692821/902500,  1932725/662596\}, \quad \{1, 2993525/2896804,   6519845/6461764\}, \\
& \{1, 3603685/2965284,   5791045/777924\}, \quad \{1, 4324625/1478656,   4919681/883600\}, \\
& \{1, 12376325/12096484,  12844709/11628100\}, \\
& \{1, 19193525/18887716,  22980245/15100996\}, \\
& \{1, 12231605/2353156,  13689845/894916\}.
\end{eqnarray*}
\end{example}

Example \ref{ex:1} suggests that there might exist infinitely many strong rational
$D(-1)$-triples containing the number $1$. We will show that this is indeed true.

\begin{example}  \label{ex:2}
Let us take a closer look at strong rational $D(-1)$-triples of the form $\{1,  689/400, c\}$.
We get the following values for $c$ with numerators and denominators less than $10^{21}$:
{\small
\begin{eqnarray*}
& 1025/64, 969425/861184, 9047825/4857616, 352915361/30030400, \\
& 109066004561/106119577600, 284429759489/271837104400, \\
& 1322025501425/1301315125504, 2253725966225/876912382096, \\
& 9055090973825/809791213456, 30776081662625/29873264334736, \\
& 41085820444721/37500436537600, 38029186636625/23706420917776, \\
& 710390547822449/245964644227600, 206973503563719329/2738904077011600, \\
& 180130335826717601/7772841524238400, 1383595259111988401/1191448219611040000, \\
& 349568886374130209/40505499045648400, 842490595967154166625/184668498086700979264.
\end{eqnarray*}
}%
\end{example}

Example \ref{ex:2} clearly indicates that we may expect that there exist
infinitely many strong rational $D(-1)$-triples of the form $\{1,  689/400, c\}$.
It is not so clear what to expect for triples of the form $\{1, 5/4, c\}$
or $\{1, 65/16, c\}$. However, as stated in Proposition \ref{ref:2},
we will show that there exist infinitely many strong rational $D(-1)$-triples of each
of these forms.

So, let $\{a,b,c\}$ be a strong rational $D(-1)$-triples with $a=1$.
Thus $b-1$, $c-1$, $b^2-1$, $c^2-1$ and $bc-1$ are perfect squares.
From the first and third condition we get $b-1=\alpha^2$, $b+1=\beta^2$ for rationals $\alpha,\beta$.
By taking $\beta^2 - 2 = \alpha^2 = (\beta - 2u)^2$, we get $\beta = \frac{2u^2+1}{2u}$ and
\begin{equation} \label{eq:b}
b= \frac{4u^4+1}{4u^2}
\end{equation}
for a non-zero rational $u$.
(If $a\neq 1$, instead of genus $0$ curve $\alpha^2+2=\beta^2$, we would have a genus $1$ curve
$\alpha^4 + 2\alpha^2 + 1-a^2=\gamma^2$.)
Analogously we get
\begin{equation} \label{eq:c}
c=  \frac{4v^4+1}{4v^2}
\end{equation}
for a non-zero rational $v$.

The only remaining condition is that $bc-1$ should be a perfect square.
By inserting (\ref{eq:b}) and (\ref{eq:c}) in $bc-1=\Box$, we get
\begin{equation} \label{eq:ec1}
(16u^4+4)v^4-16u^2v^2+4u^4+1 = z^2.
\end{equation}
This curve is a quartic in $v$ with a rational point $[u,4u^4-1]$.
Thus it can be in the standard way transformed into an elliptic curve:
\begin{equation} \label{eq:ec2}
Y^2 = X(X+32u^4+8)(X+16u^4-16u^2+4).
\end{equation}
There is a point
$$ P=[-4(4u^4+1), 16u(4u^4+1)] $$
on (\ref{eq:ec2}), which comes from the point $[u,-4u^4+1]$ on (\ref{eq:ec1}).
For all non-zero rationals $u$, the point $P$ is of infinite order
on the specialized curve (\ref{eq:ec2}) over $\mathbb{Q}$
(by Mazur's theorem \cite{Mazur}, if suffices to check that $mP \neq \mathcal{O}$ for $m\leq 12$).
Now we consider multiples $mP$, $m\geq 2$, of $P$ on (\ref{eq:ec2}),
transfer them back to the quartic (\ref{eq:ec1}),
and compute the components $b,c$ of the corresponding strong rational $D(-1)$-triple.
Since the point $P$ is of infinite order, for fixed $u$, i.e. fixed $b$,
in that way we get infinitely many strong rational $D(-1)$-triples of the form $\{1,b,c\}$,
thus proving Proposition \ref{ref:2} and Theorem \ref{ref:1}.

The point $P$ gives $[u, -4u^4+1]$, and thus does not provide a triple, since in this case we get
$v=u$ and $b=c$.
The point $2P$ gives $[-u(4u^4-3)/(12u^4-1), (64u^{12}+272u^8-68u^4-1)/(12u^4-1)^2]$
and the strong rational $D(-1)$-triple
\begin{eqnarray*}
&  \{1, (4u^4+1)/(4u^2), \\
&  (4u^4+1)(256u^{16}+4352u^{12}-1952u^8+272u^4+1)/(4u^2(4u^4-3)^2(12u^4-1)^2)
\},
\end{eqnarray*}
while $3P$ gives $[u(64u^{12}-656u^8+108u^4+5)/(320u^{12}+432u^8-164u^4+1), -(16384u^{28}+741376u^{24}-760832u^{20}+812288u^{16}-203072u^{12}+11888u^8-724u^4-1)/(320u^{12}+432u^8-164u^4+1)^2)]$
and the strong rational $D(-1)$-triple
\begin{eqnarray*}
&  \{1, (4u^4+1)/(4u^2), \\
&  (4u^4+1)(256u^{16}+4352u^{12}-195u^8+272u^4+1) \times \\  & (65536u^{32}+6422528u^{28}-13516800u^{24}+49995776u^{20}-23443968u^{16} \\
& \mbox{} +3124736u^{12}-52800u^8+1568u^4+1) \\
&  /(4u^2(64u^{12}-656u^8+108u^4+5)^2(320u^{12}+432u^8-164u^4+1)^2)
\}.
\end{eqnarray*}
By inserting $u=1$, we get the triples
$$ \{1, 5/4, 14645/484\} \quad \mbox{and} \quad \{1, 5/4, 330926870165/318391604644\} $$
(the first triple already appeared in Example \ref{ex:1}, while
the second triple is outside of the range covered by Example \ref{ex:1}).

It is clear that further multiples $4P,5P,\ldots$ would provide more complicated formulas for triples.
To get new relatively simple formulas for triples, we may try to find subfamilies of
the elliptic curve (\ref{eq:ec2}) with rank $\geq 2$.
For that purpose, we may use the method explained e.g. in \cite{Duje-Peral}.

We search for an additional point on the $2$-isogenous curve
\begin{equation} \label{eq:ec3}
y^2=x(x^2-24x+32xu^2-96xu^4+16+128u^2+384u^4+512u^6+256u^8)
\end{equation}
by considering divisors of $16 + 128u^2 + 384u^4 + 512u^6 + 256u^8=16(2u^2+1)^4$.
Imposing $x=8(2u^2+1)$ to be the $x$-coordinate of a point on (\ref{eq:ec3})
leads to the condition that $4u^2-14$ is a square, which gives
$u=(14+w^2)/(4w)$ for a rational $w$. Thus
$$ b=(w^8+56w^6+1240w^4+10976w^2+38416)/(16w^2(14+w^2)^2). $$
By transferring the additional point of infinite order of the original quartic (\ref{eq:ec1}),
we get
$$ v=(w^6+18w^4-100w^2-392)/(4w(3w^4+28w^2+140)) $$
and
\begin{eqnarray*}
&  c=(w^8+40w^6+4888w^4+7840w^2+38416) \times \\ &  (w^8-4w^7+24w^6-40w^5+152w^4+16w^3+608w^2+672w+784) \times \\ &  (w^8+4w^7+24w^6+40w^5+152w^4-16w^3+608w^2-672w+784) \\
&  /(16w^2(w^6+18w^4-100w^2-392)^2(3w^4+28w^2+140)^2)
.
\end{eqnarray*}
By inserting $w=1$, we get the triple
$\{1, 50689/3600, 104776974625/104672955024\}$
(this triple is outside of the range covered by Example \ref{ex:1}).

\section{The generic rank and generators of the families of elliptic curves} \label{sec:3}

In Section \ref{sec:2} we used families of elliptic curves with rank $\geq 1$ over
$\mathbb{Q}(u)$, resp. rank $\geq 2$ over $\mathbb{Q}(w)$, and known independent points
of infinite order to construct families of strong rational $D(-1)$-triples.
It is natural to ask what is the exact generic rank of these two families
and whether the known points are in fact generators of the corresponding Mordell-Weil groups.
The recent algorithm of Gusi\'c and Tadi\'c from \cite{GT2} (see also \cite{GT1,Stoll} for
other variants of the algorithm, and \cite{DPT,DGT} for several applications of the algorithm)
allows us to answers these questions.


First we prove that the elliptic curve given in \eqref{eq:ec2} has rank one over $\mathbb Q(u)$
and the free generator is the point $P=[-4(4u^4+1), 16u(4u^4+1)]$.
By the algorithm from \cite{GT2} we have:
\begin{itemize}
\item The specialization  at $u_0=6$ is injective by \cite[Theorem 1.1]{GT2}.

\item The coefficients  of the specialized curve are $[0,61644,0,836402720,0]$.

\item By {\tt mwrank} \cite{Cr}, the specialized curve has rank equal to $1$ and its free generator is the point $G_1=[-20740,497760]$
\item We have that the specialization of the point $P$ at $u_0=6$ satisfies
$P(6)=G_1$.
\end{itemize}
Now it is obvious that  the elliptic curve has rank $1$ and the point  $P$ is the free generator
of the elliptic curve \eqref{eq:ec2} over $\mathbb Q(u)$.

Now we consider the elliptic curve obtained from (\ref{eq:ec2}) by the substitution
$u=(14+w^2)/(4w)$. After removing the denominators, we get the elliptic curve
$C$ over $\mathbb Q(w)$ given by the equation
\begin{eqnarray*}
& & Y^2=X^3+(3 w^8+152 w^6+3272 w^4+29792 w^2+115248) X^2 \\
& & \mbox{}+2 (w^8+56w^6+1240w^4+10976w^2+38416)(w^4+20w^2+196)^2X.
\end{eqnarray*}
We claim that $C$ has  rank   equal to $2$ over $\mathbb{Q}(w)$ and that
the points with first coordinates
\begin{eqnarray*}
x(P) &\!\!=\!\!& -(w^8+56w^6+1240w^4+10976w^2+38416), \\
x(Q) &\!\!=\!\!& (w^2-14)^2(w^4+20w^2+196)^2/(64w^2)
\end{eqnarray*}
are its free generators. We again apply the algorithm from \cite{GT2}:

\begin{itemize}
\item We use the specialization  at $w_0=6$ which is injective by \cite[Theorem 1.1]{GT2}.
\item  The specialized curve over $\mathbb Q$ is $[0,17558832,0,61973480694272,0]$.

\item By {\tt mwrank} \cite{Cr}, the rank of this specialized curve over $\mathbb Q$  is
equal to $2$ and its free generators are
$$ G_1 = [2880000, 18655065600], \quad G_2 = [37002889/36, 1971840224123/216]. $$

\item
We have that for the specialization of the points $P,Q$ at $w_0=6$ it holds
$P(6) = G_1 + T$, $Q(6) = G_2$, where $T$ is a torsion point on the specialized curve.
\end{itemize}
Thus we get that the elliptic curve $C$ has rank $2$ and that the points $P$ and $Q$ are
free generators of $C$ over $\mathbb Q(w)$.

\section{Concluding remarks}

\begin{remark} \label{rem1}
We may ask how large can the rank be over $\mathbb{Q}$ of a specialization for $u\in \mathbb{Q}$
of the elliptic curve (\ref{eq:ec2}). Since for $u=(14+w^2)/(4w)$ the rank over $\mathbb{Q}(w)$
is equal to $2$, by Silverman's specialization theorem \cite[Theorem 11.4]{Sil},
we conclude that there are
infinitely many rationals $u$ for which the rank of (\ref{eq:ec2}) is $\geq 2$.
By using standard methods for searching for curves of relatively large rank
in parametric families of elliptic curves (see e.g. \cite{Duje-Glasnik}), we are able to
find curves with rank equal to $3$ (e.g. for $u=2/5$, $u=4$),
$4$ (e.g. for $u=50/11$, $u=12/65$), $5$ (e.g. for $u=12/65$, $u=16/83$) and
 $6$ (for $u=86/743$, $u=3570/1051$, $u=1642/3539$). Note that $u=2/5$ gives $b=689/400$.
 The fact that for this specialization the specialized curve has rank $3$,
 with generators with relatively small height,
 explains the observation from Example \ref{ex:2} that there are unusually many strong
 rational $D(-1)$-triples of the form $\{1, 689/400, c\}$ for $c$'s with small numerators
 and denominators (see the arguments given in \cite[Section 4]{ADKT}).
\end{remark}

\begin{remark} \label{rem2}
The results of this paper motivate following open questions:
\begin{itemize}
\item[1)] Are there infinitely many strong rational $D(-1)$-triples that do not contain the number $1$?
\item[2)] Is there any strong rational $D(-1)$-quadruple?
\end{itemize}
Note that the triple $(a,b,c)=(125/117, 689/400, 14353373/13130325)$ from Example \ref{ex:1} has an
additional property that $b-1$ is also a square. Furthermore, $26(a-1)$ and $26(c-1)$ are also squares.
Hence, although we do not know any strong $D(-1)$-quadruple over $\mathbb{Q}$, we get the set
$\{1, 125/117, 689/400, 14353373/13130325\}$ which is a strong $D(-1)$-quadruple over
the quadratic field $\mathbb{Q}(\sqrt{26})$.

\end{remark}

\begin{remark} \label{rem3}
Let $a\neq \pm 1$ be a rational such that $a^2-1$ is a square, i.e. $a=(t^2+1)/(2t)$ for a rational
$t\neq 0,\pm 1$.
It can be extended to infinitely many strong rational $D(-1)$-pairs. Indeed, as we have already mentioned,
by following the construction in the case $a=1$, we now get the condition
$\alpha^4 + 2\alpha^2 + 1-a^2=\gamma^2$. This quartic is birationally equivalent to the elliptic curve
\begin{equation} \label{eq:ect}
Y^2 = (X+2t^2)(X^2 + t^6-2t^4+t^2),
\end{equation}
for which we can show,
by taking the specialization $t_0=11$ in \cite[Theorem 1.3]{GT2}, that it has the rank over
$\mathbb{Q}(t)$ equal to $1$, with the free generator $R=[-t^2+1, t^4-1]$
(and by Mazur's theorem \cite{Mazur}, we find that $R$ is of infinite order for all rationals $t\neq 0,\pm 1$).
One explicit extension $\{a,b\}$ is by $b=(t^4+18t^2+1)/(8t(t^2+1))$.
We have noted that the elements of the known examples of strong rational $D(-1)$-triples
that do not contain the number $1$ induce the elliptic curves with relatively large rank.
In particular, for $a=42601/11849$ and $a=14353373/13130325$ the rank is equal to $5$.
We have performed an additional search for strong rational $D(-1)$-triples
that do not contain the number $1$ by considering elliptic curves in the family (\ref{eq:ect})
with the rank $\geq 3$, and checking small linear combinations of their generators. In that way,
we found one new example of the strong rational $D(-1)$-triple (corresponding to $t=17/481$):
$$ \{115825/8177, 408988121/327645760,   752442457/720825305\}, $$
which is outside of the range covered by Example \ref{ex:1}.
\end{remark}

\bigskip

{\bf Acknowledgements.}
The authors were supported by the Croatian Science Foundation under the project no.~6422. A.~D. acknowledges support from the QuantiXLie Centre of Excellence, a project
co-financed by the Croatian Government and European Union through the
European Regional Development Fund - the Competitiveness and Cohesion
Operational Programme (Grant KK.01.1.1.01.0004).

\end{document}